\documentclass[10pt]{article}
\usepackage{amssymb}
\usepackage{amsmath,amssymb}\usepackage{cite}\usepackage{mathrsfs}\usepackage[amsmath,thmmarks]{ntheorem}
\usepackage{listings}
\usepackage[titletoc]{appendix}\allowdisplaybreaks

\makeatletter
\renewcommand\theequation{\thesection.\@arabic\c@equation}
\makeatother

\newtheorem{thm}{ Theorem}[section]%
\newtheorem{lem}[thm]{ Lemma}%
\newtheorem{Exam}[thm]{ Example}%
\newtheorem{cor}[thm]{ Corollary}%

\newtheorem{Remark}[thm]{ Remark}%
\newtheorem{Pro}[thm]{ Proposition}%
\newtheorem{Def}[thm]{Definition}%
\newtheorem{Con}[thm]{Conjecture}%
\input cyracc.def

\def\f{\noindent}

\def\demo{\f{\bf Proof.}\hskip10pt}

\def\qed{\hfill $\Box$}

\begin{document}
\title{\textbf{On simultaneous conjugation of permutations}}
\author{Junyao Pan\,
 \\\\
School of Sciences, University of wuxi, Wuxi, Jiangsu,\\ 214105
 People's Republic of China \\}
\date {} \maketitle

\baselineskip=16pt

\vskip0.5cm

{\bf Abstract:} Let $\alpha$ and $\beta$ be two permutations in $S_n$. We prove that if the commutator $[\alpha,\beta]$ has at least $n-4$ fixed points then there exists a permutation $\gamma\in S_n$ such that $\alpha^\gamma=\alpha^{-1}$ and $\beta^\gamma=\beta^{-1}$. This gives an affirmative answer to a conjecture proposed by Danny Neftin, which leads to the completion of the classification of monodromy groups for sufficiently large degree.



\section {Introduction}

In this section, denote the complex by $\rm{C}$, and the projective line by $\rm{P^1}$.

Let $X,Y$ be Riemann Surfaces where $X$ has genus $g$, and $f:~X\rightarrow Y$ be a nonconstant rational map of degree $n$. The monodromy group $G$ of this cover, that is, the Galois group of the Galois closure of $\rm{C(X)/C(Y)}$. So $G$ is a transitive group of degree $n$. See \cite[~Page~1,~Chapter~1]{GS}, the main problem some scholars are interested in is the following: if we fix the genus of $X$, what restrictions are placed on $G$? For more background on this problem, refer to \cite{GN,GT}. Actually, the critical case to investigate is when $Y=\rm{P^1}$ and $f$ is indecomposable because $G$ is a primitive permutation group of degree $n$ in this case. The natural question is for which primitive group $G$ appear as the monodromy group of an indecomposable low genus covering $f:~X\rightarrow\rm{P^1}$, and for $G$ which is not alternating or symmetric (in their natural action) what are the corresponding covers $f$? This is the problem of the classification of monodromy groups, which was started by Guralnick and Thompson after the classification of finite simple groups and its roots lie in the work of Ritt and Zariski. The classification of monodromy groups is a long project with many people involved, for some results and history see \cite{GS}.

Let $B\subset\rm{P^1}$ be the set of branch points of the cover $f$ with $|B|=r$, where $f:X\rightarrow \rm{P^1}$ is indecomposable. Since the fundamental group of $\rm{P^1}\setminus B$ is a group generated by $r$ elements with the single relation that the product is $1$, using Riemann's Existence Theorem, it follows that the existence of a cover gives rise to elements $x_1,...,x_r$ in $G$ such that $G=\langle x_1,...,x_r\rangle$ and $x_1\cdot\cdot\cdot x_r=1$. Therefore, the investigation of monodromy groups associated with group structures, especially the ``primitive" Hurwitz problem of the classification of monodromy groups. In group theoretic language, it is the following problem: Given conjugacy classes $C_1,...,C_r$ in a permutation group $G$, are there elements $x_i$ in $C_i$ with product $x_1\cdot\cdot\cdot x_r=1$ and such that the group generated by $x_1,...,x_r$ is primitive, and each such tuple corresponds to an indecomposable covering which appears in this classification. Recently, the classification was completed when the degree is assumed to be sufficiently large in \cite{NZ}, except for these two families of ramification types: Given by conjugacy classes of $S_d \wr S_2$: $s, s, s, ([d-4, 2,2], 1)s$ and $s, s, s, ([d-3, 3], 1)s$, here $s$ denotes the swap in $S_d \wr S_2$, and a conjugacy classes in this group is denoted by $(a,b)s$ where $a,b$ are partition of $d$ (conjugacy classes in $S_d$). Indeed, Danny Neftin and Michael E. Zieve guessed that these two families of ramification types do not occur, and further Danny Neftin observed that if the following conjecture is true then these two families of ramification types do not occur.

\begin{Con}\label{pan1-0}\normalfont(\cite{N}~)
Let $\alpha$ and $\beta$ be two permutations in $S_n$. If the commutator $[\alpha,\beta]$ has at least $n-4$ fixed points, then there exists a permutation $\gamma$ in $S_n$ such that $\alpha^\gamma=\alpha^{-1}$ and $\beta^\gamma=\beta^{-1}$. Here $\alpha^\gamma=\gamma^{-1}\alpha\gamma$ and $[\alpha,\beta]=\alpha^{-1}\beta^{-1}\alpha\beta$.
\end{Con}

In this paper, we give a positive solution of Conjecture\ \ref{pan1-0}, and then we start to explain the idea of the proof. Firstly, we will provide some ways to construct some permutations conjugate a permutation onto its inverse under some conditions (see Lemma\ \ref{pan2-7} and Corollary\ \ref{pan2-112}) as well as a property of simultaneous conjugation (see Lemma\ \ref{pan2-8}), and then by using these ways and the property to give a formal proof of Conjecture\ \ref{pan1-0} in the case when $\alpha\beta=\beta\alpha$ (see Lemma\ \ref{pan2-6}), and further the problem will be reduced to the simpler situation (see Assumption 2.5). Secondly, we take $\alpha^{-1}$ and $\alpha^\beta$ as the objects because of $[\alpha,\beta]=\alpha^{-1}\alpha^\beta$, and we discover that there exist two situations of $\alpha^\beta$, and there exists two close relationships between the two cases, such as:

 {\bf{Example A:}} Given $\alpha_1=(1,2)(3,4)(5,6)(7,8,9)$ and $\beta_1=(1,4,6,8)(2,3,5,7)$. Then $[\alpha_1,\beta_1]=(2,9,7)$, $(1,2)^{\beta_1}=(3,4)$, $(3,4)^{\beta_1}=(5,6)$, $(5,6)^{\beta_1}=(7,8)$ and further $(2,3,5,7)^{\alpha_1}=(1,4,6,8)$.

{\bf{Example B:}} Given $\alpha_2=(1,2)(7,8,9)$, $\beta_2=(1,8)(2,7)$. Then $[\alpha_2,\beta_2]=[\alpha_1,\beta_1]$, $(1,2)^{\beta_2}\neq(1,2)$, $(1,2)^{\beta_2}\neq(7,8,9)$, $(7,8,9)^{\beta_2}\neq(1,2)$ and $(7,8,9)^{\beta_2}\neq(7,8,9)$.

{\bf{Relationship 1:}} There exists a permutation $\mu=(6,4,1)(5,3,2)$ such that $\alpha_1^\mu=\alpha_1$ and $\beta_2=\mu\beta_1$ and $3,4,5,6$ are fixed points of $\beta_2$.

{\bf{Relationship 2:}} Pick $\omega=(1,2)(7,8)$ and $\nu=(4,5)(3,6)$ such that $\alpha_2^{\omega}=\alpha_2^{-1}$, $(1,2)^\omega=(1,2)=(1,2)^{-1}$, $(7,8,9)^\omega=(8,7,9)=(7,8,9)^{-1}$ and $\beta_2^{\omega}=\beta_2^{-1}$ and $\mu^{\omega\nu}=(3,5,2)(4,6,1)=\mu^{-1}$. Then we see that $\alpha_1^{\omega\nu\mu}=(6,5)(4,3)(1,2)(8,7,9)=\alpha_1^{-1}$ and $\beta_1^{\omega\nu\mu}=(5,3,2,7)(6,4,1,8)=\beta_1^{-1}$.\\
So we will introduce a definition (See Definition\ \ref{pan2-000000}) to distinguish between the two situations, and the two relationships will be summarized as three lemmas (see Lemma\ \ref{pan4-004}, Lemma\ \ref{pan4-0004} and Lemma\ \ref{pan4-00004}) which imply that the case of Example A can be transformed into the case of Example B, and so we shall introduce a tool (See Definition\ \ref{pan2-0}) to characterize the case of Example B (see Theorems\ \ref{pan3-6}~-\ \ref{pan3-20}). Finally, we wll apply the techniques (see Lemma\ \ref{pan2-66} and Lemma\ \ref{pan2-12}) to confirm Conjecture\ \ref{pan1-0} for the case of Example B (see Lemma\ \ref{pan4-00400}), and then by using lemmas (see Lemma\ \ref{pan4-00004} and Lemma\ \ref{pan4-00400}) to confirm Conjecture\ \ref{pan1-0} for the case of Example A (see Lemma\ \ref{pan4-0}).

\section {Preliminaries}

Throughout this paper, if there is no special statement, the different letters and the same letter with different marks indicate different points. Furthermore, $S_n$ is the symmetric group on the set $[n]=\{1,2,...,n\}$, and $(1)$ is the identity of $S_n$.

Recall some notions and notations about symmetric group which will be used, for further details see\cite{C,D}. For $\alpha,\beta\in S_n$, we write $Fix(\alpha)=\{i\in[n]|i^\alpha=i\}$ and $M(\alpha)=[n]\setminus Fix(\alpha)$ and the commutator $[\alpha,\beta]=\alpha^{-1}\beta^{-1}\alpha\beta=\alpha^{-1}\alpha^\beta$ and $o(\alpha)$ the order of $\alpha$. It is well known that there exist some disjoint cycles $\alpha_1,...,\alpha_s$ such that $\alpha=\alpha_1\cdot\cdot\cdot\alpha_s$, and then we use $\{\alpha\}$ to denote the set $\{\alpha_1,...,\alpha_s\}$, and we say $\alpha_i$ is a cycle factor of $\alpha$. It is obvious that $\alpha^{-1}=\alpha^{-1}_1\cdot\cdot\cdot\alpha^{-1}_s$ and $\alpha^\beta=\alpha^{\beta}_1\cdot\cdot\cdot\alpha^{\beta}_s$ and $\{\alpha^\beta\}=\{\alpha^{\beta}_1,...,\alpha^{\beta}_s\}$. Moreover, the centralizer of $\alpha$ in $S_n$ is denoted by $C_{S_n}(\alpha)$, that is, $C_{S_n}(\alpha)=\{x\in S_n|x\alpha=\alpha x\}$. In particular, we note that for $\gamma,\omega\in S_n$, $\alpha^\gamma=\alpha^\omega$ if and only if $\gamma\omega^{-1}\in C_{S_n}(\alpha)$. So we start with the $\omega$ such that $\alpha_i^\omega=\alpha_i^{-1}$ for each $i=1,2,...,s$.

\begin{lem}\label{pan2-7}\normalfont
Let $\alpha=(x_r,x_{r-1},...,x_1)$ and $\beta=(y_r,y_{r-1},...,y_1)$ be two disjoint $r$-cycles in $S_n$.

(i) For each $x_i$ in $M(\alpha)$, there exists an involution $\omega$ such that $\alpha^{\omega}=\alpha^{-1}$ and $x_i\in Fix(\omega)$ and $M(\omega)\subseteq M(\alpha)$. In particular, there are at most two points in $M(\alpha)\cap Fix(\omega)$.

(ii) For any two points $x_i$ and $x_j$ in $M(\alpha)$, there exists an involution $\omega$ such that $\alpha^{\omega}=\alpha^{-1}$ and $x_i^\omega=x_j$ and $M(\omega)\subseteq M(\alpha)$.

(iii) For all $x_i\in M(\alpha)$ and $y_j\in M(\beta)$, there exists an involution $\omega$ such that $\alpha^\omega=\beta^{-1}$, $\beta^{\omega}=\alpha^{-1}$, $M(\omega)\subseteq M(\alpha)\cup M(\beta)$ and $(x_i,y_j)\in\{\omega\}$.

\end{lem}
\demo Note that $\alpha^{-1}=(x_1,x_2,...,x_r)$ and $\alpha=(x_k,x_{k-1},...,x_1,x_r,x_{r-1},...,x_{k+1})$ and  $\alpha^{\omega}=(x_k^\omega,x_{k-1}^\omega,...,x_1^\omega,x_r^\omega,x_{r-1}^\omega,...,x_{k+1}^\omega)$ where $1\leq k\leq r$. So we see that $\alpha^{\omega}=\alpha^{-1}$ if and only if there exists an integer $k$ such that $x_k^\omega=x_1,x_{k-1}^\omega=x_{2},...,x_1^\omega=x_k,x_{r}^\omega=x_{k+1},...,x_{k+1}^\omega=x_{r}$. Picking an involution $\omega$ such that $\{\omega\}=\{(x_s,x_t)|x_s,x_t\in M(\alpha)~{\rm{with}}~s+t=k+1~{\rm{or}}~k+r+1\}$. One easily checks that $\alpha^{\omega}=\alpha^{-1}$ and $M(\omega)\subseteq M(\alpha)$, and further the $\omega$ is completely determined by $k$.

(i) According to the above discussions, it follows that $x_i\in M(\alpha)\cap Fix(\omega)$ if and only if there exists an integer $k$ such that $1\leq k\leq r$, $2i=k+1$ or $2i=k+r+1$. Since $1\leq i\leq r$, there always exists such integer $k$ for each $x_i$ in $M(\alpha)$. In addition, the necessary and sufficient condition shows that there are at most two points in $M(\alpha)\cap Fix(\omega)$, and $|M(\alpha)\cap Fix(\omega)|=2$ if and only if both $k+1$ and $k+r+1$ are even numbers.

(ii) Similarly, we know that $(x_i,x_j)\in\{\omega\}$ if and only if there exists an integer $k$ such that $1\leq k\leq r$, $i+j=k+1$ or $i+j=k+r+1$. It follows from $1\leq i,j\leq r$ that there always exists such integer $k$ for any two points $x_i$ and $x_j$.

(iii) Note that $\alpha=(x_i,x_{i-1},...,x_1,x_r,x_{r-1},...,x_{i+1})$ and $\beta^{-1}=(y_j,y_{j+1},...,y_r,y_1,...,y_{j-1})$. Since $\alpha$ and $\beta$ are two disjoint $r$-cycles, we can take the involution $\omega$ such that $x_i^\omega=y_j,x_{i-1}^\omega=y_{j+1},...x_{i+1}^\omega=y_{j-1}$ and $M(\omega)= M(\alpha)\cup M(\beta)$. It is simple to see $\alpha^\omega=\beta^{-1}$ and $(x_i,y_j)\in\{\omega\}$. In addition, we see $\alpha=\alpha^{\omega^2}=(\beta^{-1})^\omega$, and so $\beta^{\omega}=\alpha^{-1}$. We have thus proved this lemma. \qed

\begin{cor}\label{pan2-112}\normalfont
Let $\alpha=(x_{11},x_{12},...,x_{1r})\cdot\cdot\cdot(x_{k1},x_{k2},...,x_{kr})$ and $\beta=(x_{11},x_{21},...,x_{k1})$. Then for any permutation $\omega$ with $\beta^{\omega}=\beta^{-1}$ and $M(\omega)\subseteq M(\beta)$, there exists a permutation $\nu$ such that $\alpha^{\omega\nu}=\alpha^{-1}$ and $M(\nu)\subseteq M(\alpha)\setminus M(\beta)$ and $\beta^{\omega\nu}=\beta^{-1}$.

\end{cor}
\demo Let $\alpha_i=(x_{i1},x_{i2},...,x_{ir})\in\{\alpha\}$ for $i=1,2,...,k$. From $\beta^{\omega}=\beta^{-1}$ and $M(\omega)\subseteq M(\beta)$, it follows that $\omega$ is an involution, and so we assume that $$\omega=(x_{s_11},x_{t_11})(x_{s_21},x_{t_21})\cdot\cdot\cdot(x_{s_l1},x_{t_l1})~{\rm{where}}~\{\omega\}=\{(x_{s_11},x_{t_11}),(x_{s_21},x_{t_21}),...,(x_{s_l1},x_{t_l1})\}.$$ For each cycle factor $(x_{s_j1},x_{t_j1})$, by Lemma\ \ref{pan2-7} (iii), there exists an involution $\omega_j$ such that $$\alpha_{s_j}^{(x_{s_j1},x_{t_j1})\omega_j}=\alpha_{t_j}^{-1}~{\rm{and}}~\alpha_{t_j}^{(x_{s_j1},x_{t_j1})\omega_j}=\alpha_{s_j}^{-1}~{\rm{and}}~M(\omega_j)\subseteq M(\alpha_{s_j}\alpha_{t_j})\setminus\{x_{s_j1},x_{t_j1}\}.$$ Picking $\nu'=\omega_1\omega_2\cdot\cdot\cdot\omega_l$. Then we see that $$M(\nu')\subseteq M(\alpha)\setminus M(\beta)~{\rm{and}}~(\alpha_{{s_1}}\alpha_{{s_2}}\cdot\cdot\cdot\alpha_{{s_l}}\alpha_{{t_1}}\alpha_{{t_2}}\cdot\cdot\cdot\alpha_{{t_l}})^{\omega\nu'}=(\alpha_{{s_1}}\alpha_{{s_2}}\cdot\cdot\cdot\alpha_{{s_l}}\alpha_{{t_1}}\alpha_{{t_2}}\cdot\cdot\cdot\alpha_{{t_l}})^{-1}.$$

On the other hand, Lemma\ \ref{pan2-7} (i) shows that there are at most two points in $M(\beta)\cap Fix(\omega)$. If $M(\beta)\cap Fix(\omega)=\emptyset$, then $\alpha=\alpha_{{s_1}}\alpha_{{s_2}}\cdot\cdot\cdot\alpha_{{s_l}}\alpha_{{t_1}}\alpha_{{t_2}}\cdot\cdot\cdot\alpha_{{t_l}}$, and so $\alpha^{\omega\nu}=\alpha^{-1}$ and $M(\nu)\subseteq M(\alpha)\setminus M(\beta)$ and $\beta^{\omega\nu}=\beta^{-1}$ where $\nu=\nu'$. Suppose $\{x_{p1}\}= M(\beta)\cap Fix(\omega)$. Then by Lemma\ \ref{pan2-7} (i) we see that there exists an involution $\omega'$ such that $\alpha_p^{\omega'}=\alpha_p^{-1}$ and $M(\omega')\subseteq M(\alpha_p)\setminus\{x_{p1}\}$. In this case, $\alpha=(\alpha_{{s_1}}\alpha_{{s_2}}\cdot\cdot\cdot\alpha_{{s_l}}\alpha_{{t_1}}\alpha_{{t_2}}\cdot\cdot\cdot\alpha_{{t_l}})\alpha_p$, and it is easy to verify $\alpha^{\omega\nu}=\alpha^{-1}$ and $M(\nu)\subseteq M(\alpha)\setminus M(\beta)$ and $\beta^{\omega\nu}=\beta^{-1}$ where $\nu=\omega'\nu'$. Assume $\{x_{p1},x_{q1}\}=M(\beta)\cap Fix(\omega)$. Similarly, there also exists an involution $\omega''$ such that $\alpha_q^{\omega''}=\alpha_q^{-1}$ and $M(\omega'')\subseteq M(\alpha_q)\setminus\{x_{q1}\}$, taking $\nu=\omega'\omega''\nu'$, we see $\alpha^{\omega\nu}=\alpha^{-1}$ and $M(\nu)\subseteq M(\alpha)\setminus M(\beta)$ and $\beta^{\omega\nu}=\beta^{-1}$. The proof of this corollary is complete. \qed

Next, we first state an useful lemma and then give a formal proof of Conjecture\ \ref{pan1-0} in the case when $\alpha\beta=\beta\alpha$.

\begin{lem}\label{pan2-8}\normalfont
Let $\alpha$ and $\beta$ be two permutations in $S_n$. If there exists a permutation $\mu\in S_n$ such that $\alpha^\mu=\alpha^{-1}$ and $\beta^\mu=\beta^{-1}$, then there exists a permutation $\gamma\in S_n$ such that $\alpha^\gamma=\alpha^{-1}$, $\beta^\gamma=\beta^{-1}$ and $M(\gamma)\subseteq M(\alpha)\cup M(\beta)$.
\end{lem}
\demo We define a subset of $\{\mu\}$, that is, $IN(\alpha,\beta)=\{\mu_i\in\{\mu\}|~M(\mu_i)\subseteq M(\alpha)\cup M(\beta)\}$. If $IN(\alpha,\beta)\neq\emptyset$, it is easy to see that $\alpha^\gamma=\alpha^{-1}$, $\beta^\gamma=\beta^{-1}$ and $M(\gamma)\subseteq M(\alpha)\cup M(\beta)$, where $\{\gamma\}=IN(\alpha,\beta)$. If $IN(\alpha,\beta)=\emptyset$, then $M(\mu)\cap M(\alpha)=\emptyset$ and $M(\mu)\cap M(\beta)=\emptyset$, and thus $\alpha=\alpha^{-1}$ and $\beta=\beta^{-1}$, which imply $\alpha$ and $\beta$ are two involutions. In this case, we may set $\gamma=(1)$ where $(1)$ is the identity. It is clear that $\alpha^\gamma=\alpha^{-1}$, $\beta^\gamma=\beta^{-1}$ and $M(\gamma)=\emptyset\subseteq M(\alpha)\cup M(\beta)$. The proof of this lemma is now complete.   \qed

\begin{lem}\label{pan2-6}\normalfont({\cite[Danny Neftin]{N}}~)
Let $\alpha$ and $\beta$ be two permutations in $S_n$ with $\alpha\beta=\beta\alpha$. Then there exists a permutation $\gamma\in S_n$ such that $\alpha^\gamma=\alpha^{-1}$ and $\beta^\gamma=\beta^{-1}$.
\end{lem}
\demo Since $\alpha\beta=\beta\alpha$, it follows that $\{\alpha\}$ is an invariant set under the action of $\beta$ by
conjugation. According to Lemma\ \ref{pan2-8}, it suffices to prove the case that $\{\alpha\}$ is an orbit under the action of $\beta$ by
conjugation and vice versa. If $\alpha$ is a cycle, then $\beta=\alpha^{m}$ for a positive integer $m$, and it is simple to show that there exists a permutation which simultaneously conjugates $\alpha$ and $\beta$ onto their respective inverses. So we let $\alpha=\alpha_1\alpha_2\cdot\cdot\cdot\alpha_k$, and the cycle factor $\alpha_i=(x_{i1},x_{i2},...,x_{il})$ for $i=1,2,...,k$ with $k>1$. Since $\{\alpha\}$ is an orbit of $\beta$ which acts on it by conjugation, we may assume that $\alpha_k^\beta=\alpha_{1}$ and $\alpha_i^\beta=\alpha_{i+1}$ for $i=1,2,...,k-1$. Without loss of generality, we set $x^\beta_{ij}=x_{(i+1)j}$ for $i=1,2,...,k-1$ and $j=1,2,...,l$.
Picking $\mu=\mu_1\mu_2\cdot\cdot\cdot\mu_l$, where $\mu_j=(x_{kj},x_{(k-1)j},...,x_{1j})$ for $j=1,2,...,l$. Then we see that $\alpha_1^\mu=\alpha_{k}$, $\alpha_i^\mu=\alpha_{i-1}$ for $i=2,...,k$, $M(\mu)=M(\alpha)$ and further $M(\alpha_2\cdot\cdot\cdot\alpha_k)\subseteq Fix(\mu\beta)$. Obviously, $\alpha_1^{\mu\beta}=\alpha_1$, and thus $\mu\beta=\alpha_1^{m}$ for a positive integer $m$, therefore, there exists an involution $\omega$ so that $\alpha_1^\omega=\alpha_1^{-1}$ and $(\mu\beta)^\omega=(\mu\beta)^{-1}$ and $M(\omega)\subseteq M(\alpha_1)$. Note that $\alpha_1=(x_{11},x_{12},...,x_{1l})$ and $\mu=(x_{k1},x_{(k-1)1},...,x_{11})\cdot\cdot\cdot(x_{kl},x_{(k-1)l},...,x_{1l})$. Then by Corollary\ \ref{pan2-112} we see that there exists a permutation $\nu$ such that $\mu^{\omega\nu}=\mu^{-1}$ and $M(\nu)\subseteq M(\alpha_2\cdot\cdot\cdot\alpha_k)$. So we have $(\mu\beta)^{\omega\nu}=\mu^{\omega\nu}\beta^{\omega\nu}=\mu^{-1}\beta^{\omega\nu}$ and $(\mu\beta)^{\omega\nu}=(\alpha_1^{m})^{\omega\nu}=\alpha_1^{-m}=(\mu\beta)^{-1}=\beta^{-1}\mu^{-1}$, therefore, $\mu^{-1}\beta^{\omega\nu}=\beta^{-1}\mu^{-1}$ and so $\beta^{\omega\nu\mu}=\beta^{-1}.$ Note $\alpha=\alpha^\mu_1\alpha^{\mu^2}_1\cdot\cdot\cdot\alpha^{\mu^k}_1$. On the other hand, we see that $$(\alpha_1^{\mu^i})^{\omega\nu\mu}=\alpha_1^{\mu^i\omega\nu\mu}=\alpha_1^{(\omega\nu)(\omega\nu)^{-1}\mu^i\omega\nu\mu}=(\alpha^{-1}_1)^{(\mu^i)^{\omega\nu}\mu}=(\alpha_1^{-1})^{\mu^{-i+1}}=(\alpha_1^{\mu^{-i+1}})^{-1},$$ and thus $\alpha^{\omega\nu\mu}=(\alpha^{\mu^{-1+1}}_1\alpha^{\mu^{-2+1}}_1\cdot\cdot\cdot\alpha^{\mu^{-k+1}}_1)^{-1}$. Since $o(\mu)=k$, we have $$\alpha^{\omega\nu\mu}=(\alpha^{\mu^{k-1+1}}_1\alpha^{\mu^{k-2+1}}_1\cdot\cdot\cdot\alpha^{\mu^{k-k+1}}_1)^{-1}=(\alpha^\mu_1\alpha^{\mu^2}_1\cdot\cdot\cdot\alpha^{\mu^k}_1)^{-1}=\alpha^{-1}.$$ Consequently, we derive this lemma. \qed

Then by Lemma\ \ref{pan2-8} and Lemma\ \ref{pan2-6}, we conclude this section by the following assumption.

{\bf{Assumption 2.5:}} we assume that there is no orbit of $\beta$ acting on $\{\alpha\}$ under conjugation in the rest of this paper.

\section {Characterization of $\alpha$ and $\beta$ with $|Fix([\alpha,\beta])|\geq n-4$}

Given $\alpha=\alpha_1\cdot\cdot\cdot\alpha_s$ with $s\geq1$. Then $[\alpha,\beta]=\alpha_1^{-1}\cdot\cdot\cdot\alpha_s^{-1}\alpha_1^\beta\cdot\cdot\cdot\alpha_s^\beta$. Since $|Fix([\alpha,\beta])|\geq n-4$, it follows that there are some points which are not only in $M(\alpha)$ but also in $Fix([\alpha,\beta])$ when $|M(\alpha)|>4$. Note that if there exists an $\alpha_j$ such that $\alpha_j^\beta=\alpha_i$, then $x\in Fix([\alpha,\beta])$ for all $x\in M(\alpha_i)$. On the other hand, Assumption 2.5 implies that there exists a positive integer $k\geq2$ such that $\alpha_j^{\beta^k}\notin\{\alpha\}$. In fact, this situation really may arise, see Example A. So we introduce the following definition to characterize this situation.

\begin{Def}\label{pan2-000000}\normalfont
Let $\alpha$ and $\beta$ be two permutations in $S_n$, and $\{\alpha_1,\alpha_2,...,\alpha_k\}\subseteq\{\alpha\}$ with $k\geq2$. If $\alpha_i^{\beta}=\alpha_{i+1}$ for all $1\leq i<k$, $\alpha_k^{\beta}\notin\{\alpha\}$ and $\alpha_j^\beta\neq\alpha_1$ for each $\alpha_j\in\{\alpha\}$, then we say $\alpha_{1}\rightarrow\alpha_2\rightarrow\cdot\cdot\cdot\rightarrow\alpha_k$ is a \emph{transitive cycle factors chain} of $\beta$ on $\alpha$, and $k$ is the \emph{length} of this transitive cycle factors chain.
\end{Def}

Then by Definition\ \ref{pan2-000000} we obtain the following remark immediately.

\begin{Remark}\label{pan2-000006}\normalfont
Let $\alpha,\beta\in S_n$ with some transitive cycle factors chains of $\beta$ on $\alpha$. Then there is no common cycle factor between two different transitive cycle factors chains of $\beta$ on $\alpha$.
\end{Remark}

Inspired by Relationship 1 and Relationship 2, we get the following three lemmas.

\begin{lem}\label{pan4-004}\normalfont
Let $\alpha,\beta\in S_n$ and $\alpha_1\rightarrow\alpha_2\rightarrow\cdot\cdot\cdot\rightarrow\alpha_k$ be a transitive cycle factors chain of $\beta$ on $\alpha$. Then there exists a permutation $\mu\in C_{S_n}(\alpha_1\alpha_2\cdot\cdot\cdot\alpha_k)$ such that $o(\mu)=k$, $M(\mu)=M(\alpha_1\alpha_2\cdot\cdot\cdot\alpha_k)$ and $M(\alpha_2\alpha_3\cdot\cdot\cdot\alpha_k)\subseteq Fix(\mu\beta)$.
\end{lem}
\demo Since $\alpha_1\rightarrow\alpha_2\rightarrow\cdot\cdot\cdot\rightarrow\alpha_k$ is a transitive cycle factors chain of $\beta$ on $\alpha$, we have $\alpha^\beta_i=\alpha_{i+1}$ for $i=1,2,...,k-1$. Without loss of generality, we may assume that $\alpha_i=(x_{i1},x_{i2},...,x_{il})$ for $i=1,2,...,k$, and $x^\beta_{ij}=x_{(i+1)j}$ for $i=1,2,...,k-1$ and $j=1,2,...,l$. Pick $\mu=\mu_1\mu_2\cdot\cdot\cdot\mu_l$ where $\mu_j=(x_{kj},x_{(k-1)j},...,x_{1j})$ for $j=1,2,...,l$. Clearly, $\mu\in C_{S_n}(\alpha_1\alpha_2\cdot\cdot\cdot\alpha_k)$, $o(\mu)=k$ and $M(\mu)=M(\alpha_1\alpha_2\cdot\cdot\cdot\alpha_k)$. Moreover, we see $x^{\mu\beta}_{ij}=x_{ij}$ for $i=2,...,k$ and $j=1,2,...,l$, and thus $M(\alpha_2\alpha_3\cdot\cdot\cdot\alpha_k)\subseteq Fix(\mu\beta)$, as desired.  \qed

\begin{lem}\label{pan4-0004}\normalfont
Let $\alpha,\beta,\mu\in S_n$ with $\mu\in C_{S_n}(\alpha_1\alpha_2\cdot\cdot\cdot\alpha_k)$, $o(\mu)=k$, $M(\mu)=M(\alpha_1\alpha_2\cdot\cdot\cdot\alpha_k)$ and $M(\alpha_2\cdot\cdot\cdot\alpha_k)\subseteq Fix(\mu\beta)$, where $\alpha_1\rightarrow\alpha_2\rightarrow\cdot\cdot\cdot\rightarrow\alpha_k$ is a transitive cycle factors chain of $\beta$ on $\alpha$. Then $[\alpha,\beta]=[\alpha'\alpha_1,\mu\beta]$, where $\alpha'$ is a permutation with $\{\alpha'\}=\{\alpha\}\setminus\{\alpha_1,\alpha_2,...,\alpha_k\}$.
\end{lem}
\demo Note $[\alpha,\beta]=[\alpha'\alpha_1\alpha_2\cdot\cdot\cdot\alpha_k,\beta]=(\alpha'\alpha_1\alpha_2\cdot\cdot\cdot\alpha_k)^{-1}(\alpha'\alpha_1\alpha_2\cdot\cdot\cdot\alpha_k)^\beta$. It follows from $\alpha^\mu=\alpha$ and $M(\alpha_2\alpha_3\cdot\cdot\cdot\alpha_k)\subseteq Fix(\mu\beta)$ that $(\alpha'\alpha_1\alpha_2\cdot\cdot\cdot\alpha_k)^\beta=(\alpha'\alpha_1\alpha_2\cdot\cdot\cdot\alpha_k)^{\mu\beta}=(\alpha'\alpha_1)^{\mu\beta}\alpha_2\cdot\cdot\cdot\alpha_k$. Hence, we have $[\alpha,\beta]=[(\alpha'\alpha_1)^{-1}\alpha^{-1}_2\cdot\cdot\cdot\alpha^{-1}_k][\alpha_2\cdot\cdot\cdot\alpha_k(\alpha'\alpha_1)^{\mu\beta}]=(\alpha'\alpha_1)^{-1}(\alpha'\alpha_1)^{\mu\beta}=[\alpha'\alpha_1,\mu\beta]$. \qed

\begin{lem}\label{pan4-00004}\normalfont
Let $\alpha,\alpha',\beta,\mu\in S_n$ with $\mu\in C_{S_n}(\alpha_1\alpha_2\cdot\cdot\cdot\alpha_k)$, $o(\mu)=k$, $M(\mu)=M(\alpha_1\alpha_2\cdot\cdot\cdot\alpha_k)$ and $M(\alpha_2\cdot\cdot\cdot\alpha_k)\subseteq Fix(\mu\beta)$, where $\alpha_1\rightarrow\alpha_2\rightarrow\cdot\cdot\cdot\rightarrow\alpha_k$ is a transitive cycle factors chain of $\beta$ on $\alpha$ and $\{\alpha'\}=\{\alpha\}\setminus\{\alpha_1,\alpha_2,...,\alpha_k\}$. If there exists a permutation $\omega$ such that $M(\omega)\subseteq M(\alpha'\alpha_1)\cup M(\mu\beta)$ and $(\mu\beta)^\omega=(\mu\beta)^{-1}$ and $\alpha'^\omega=\alpha'^{-1}$ and $\alpha_1^\omega=\alpha_1^{-1}$, then there exists a permutation $\gamma$ such that $\alpha^\gamma=\alpha^{-1}$ and $\beta^\gamma=\beta^{-1}$. In particular, $\alpha'^\gamma=\alpha'^{-1}$.
\end{lem}
\demo Since $\alpha_1\rightarrow\alpha_2\rightarrow\cdot\cdot\cdot\rightarrow\alpha_k$ is a transitive cycle factors chain of $\beta$ on $\alpha$, we may set $\alpha_i=(x_{i1},x_{i2},...,x_{il})$ for $i=1,2,...,k$. As with $\mu\in C_{S_n}(\alpha_1\alpha_2\cdot\cdot\cdot\alpha_k)$, $o(\mu)=k$, $M(\mu)=M(\alpha_1\alpha_2\cdot\cdot\cdot\alpha_k)$ and $M(\alpha_2\cdot\cdot\cdot\alpha_k)\subseteq Fix(\mu\beta)$, we may take $\mu=\mu_1\mu_2\cdot\cdot\cdot\mu_l$ where $\mu_j=(x_{kj},x_{(k-1)j},...,x_{1j})$ for $j=1,2,...,l$. Note $\mu=(x_{k1},x_{(k-1)1},...,x_{11})\cdot\cdot\cdot(x_{kl},x_{(k-1)l},...,x_{1l})$ and $\alpha_1=(x_{11},x_{12},...,x_{1l})$. Then by Corollary\ \ref{pan2-112} we see that there exists a permutation $\nu$ such that $\mu^{\omega\nu}=\mu^{-1}$ and $M(\nu)\subseteq M(\alpha_2\cdot\cdot\cdot\alpha_k)$. From $M(\alpha_2\cdot\cdot\cdot\alpha_k)\subseteq Fix(\mu\beta)$, we derive $(\mu\beta)^{\omega\nu}=((\mu\beta)^{-1})^\nu=(\mu\beta)^{-1}=\beta^{-1}\mu^{-1}$. In addition, $(\mu\beta)^{\omega\nu}=\mu^{\omega\nu}\beta^{\omega\nu}=\mu^{-1}\beta^{\omega\nu}$. Thus we have $\mu^{-1}\beta^{\omega\nu}=\beta^{-1}\mu^{-1}$, which indicates $\mu^{-1}\beta^{\omega\nu}\mu=\beta^{-1}$, and so $\beta^{\omega\nu\mu}=\beta^{-1}$.

We claim that $\alpha^{\omega\nu\mu}=\alpha^{-1}$. Clearly,  $\alpha^{\omega\nu\mu}=\alpha'^{\omega\nu\mu}(\alpha_1\alpha_2\cdot\cdot\cdot\alpha_k)^{\omega\nu\mu}$. Since $M(\alpha')\cap M(\nu\mu)=\emptyset$, we have $\alpha'^{\omega\nu\mu}=\alpha'^{-1}$. Note $\alpha_1\alpha_2\cdot\cdot\cdot\alpha_k=\alpha^{\mu^1}_1\alpha^{\mu^2}_1\cdot\cdot\cdot\alpha^{\mu^{k}}_1$. Consider $(\alpha_1^{\mu^i})^{\omega\nu\mu}$ for each $i=1,2,...,k$. We see $(\alpha_1^{\mu^i})^{\omega\nu\mu}=\alpha_1^{\mu^i\omega\nu\mu}=\alpha_1^{\omega\nu(\omega\nu)^{-1}\mu^i\omega\nu\mu}$. Since $\alpha_1^{\omega}=\alpha_1^{-1}$ and $M(\nu)\cap M(\alpha_1)=\emptyset$, we have $\alpha_1^{\omega\nu(\omega\nu)^{-1}\mu^i\omega\nu\mu}=(\alpha_1^{-1})^{(\mu^i)^{\omega\nu}\mu}$. As $\mu^{\omega\nu}=\mu^{-1}$, we derive $(\alpha_1^{-1})^{(\mu^i)^{\omega\nu}\mu}=(\alpha_1^{-1})^{\mu^{-i+1}}=(\alpha_1^{\mu^{-i+1}})^{-1}$. Since $o(\mu)=k$, we have $(\alpha_1^{\mu^{-i+1}})^{-1}=(\alpha_1^{\mu^{k-i+1}})^{-1}$. Hence, $(\alpha_1\alpha_2\cdot\cdot\cdot\alpha_k)^{\omega\nu\mu}=(\alpha^{\mu^1}_1\alpha^{\mu^2}_1\cdot\cdot\cdot\alpha^{\mu^{k}}_1)^{\omega\nu\mu}=(\alpha^{\mu^1}_1\alpha^{\mu^2}_1\cdot\cdot\cdot\alpha^{\mu^{k}}_1)^{-1}=(\alpha_1\alpha_2\cdot\cdot\cdot\alpha_k)^{-1}$.
In particular, $M(\nu)\cap M(\alpha')=\emptyset$ and $M(\mu)\cap M(\alpha')=\emptyset$ indicate $\alpha'^{\omega\nu\mu}=\alpha'^{-1}$.   \qed

It follows from Lemma\ \ref{pan4-0004} and Lemma\ \ref{pan4-00004} that the problem can be attributed to the case that there is no transitive cycle factors chain. So we first deal with the case that there is no transitive cycle factors chain of $\beta$ on $\alpha$. See Example B, we have $\alpha_2^{-1}=(1,2)(8,7,9)$ and $\alpha_2^{\beta_2}=(8,7)(2,1,9)$ and $1,8\in Fix([\alpha_2,\beta_2])$ even if there is no transitive cycle factors chain of $\beta_2$ on $\alpha_2$. In fact, we had observed this phenomenon and further applied it in \cite{P}. Now we recall some notions and notations which will be used to characterize this phenomenon, for further details refer to \cite{MR}.

A \emph{block} $A$ in a cycle $\pi=(a_1,a_2,...,a_m)$ is a consecutive nonempty substring $a_i,a_{i+1},...,a_{i+l}$ of $a_i,a_{i+1},...,a_m,a_1,a_2,...,a_{i-1}$ which means that for each $1\leq k\leq l$, $a_{i+k}=a_{i+k}$ if $i+k\leq m$, and $a_{i+k}=a_j$ if $i+k>m$, where $j\equiv i+k(mod ~m)$ with $1\leq j<m$. We say two blocks $A$ and $B$ are \emph{disjoint} if they do not have points in common, and the \emph{product} $AB$ of two disjoint blocks is defined as the usual concatenation of strings, for example, if $A=1,2,3$ and $B=4,5,6,7$ then $AB=1,2,3,4,5,6,7$. A \emph{block partition} of the cycle $\pi$ is a set $\{A_1,...,A_l\}$ of pairwise disjoint blocks in $\pi$ such that there exist a block product $A_{i_1}\cdot\cdot\cdot A_{i_l}$ of these blocks such that $\pi=(A_{i_1}\cdot\cdot\cdot A_{i_l})$.

 \begin{Def}\label{pan2-0}\normalfont
 Let $g$ and $h$ be two permutations in $S_n$. If $A=x_1,x_2,...,x_l$ is a block in a cycle factor of $g$, and $A^{-1}=x_l,x_{l-1},...,x_1$ is a block in a cycle factor of $h$ such that $x_l^g\neq x_l^{h^{-1}}$ and $x_1^{g^{-1}}\neq x_1^h$, then we say $\{A,A^{-1}\}$ is a \emph{local inverse pair} between $g$ and $h$. In particular, if $l=1$, then $x_1^g\neq x_1^{h^{-1}}$ and $x_1^{g^{-1}}\neq x_1^h$. Furthermore, we say $x$ is a \emph{free point} between $g$ and $h$ if $x$ is either in $M(g)\cap Fix(h)$ or in $M(h)\cap Fix(g)$.
 \end{Def}

\begin{Remark}\label{pan3-2046}\normalfont
Let $g$ and $h$ be two permutations in $S_n$. If $|M(g)|=|M(h)|$, then the number of free points between $g$ and $h$ is even.
\end{Remark}
\demo Let $F(g,h)$ be the set of free points between $g$ and $h$. As with Definition\ \ref{pan2-0}, we have $F(g,h)=[M(g)\setminus(M(g)\cap M(h))]\cup[M(h)\setminus(M(g)\cap M(h))]$. In addition, it can easily be seen that $[M(g)\setminus(M(g)\cap M(h))]\cap[M(h)\setminus(M(g)\cap M(h))]=\emptyset$, and so $|F(g,h)|=2|M(g)\setminus(M(g)\cap M(h))|$ due to $|M(g)|=|M(h)|$. Thus we have proved this remark. \qed

\begin{lem}\label{pan2-2}\normalfont
Let $g=g_1\cdot\cdot\cdot g_s$ and $h=h_1\cdot\cdot\cdot h_t$ with $g_ih_j\neq(1)$ for all $i$ and $j$.

(i) If $x\in M(g)\cap Fix(gh)$, then there exists a local inverse pair containing $x$ between $g$ and $h$.

(ii) For each local inverse pair $\{A,A^{-1}\}$ between $g$ and $h$, there exists only one point which is not only in $A$ but also in $M(gh)$ while other points of $A$ are in $Fix(gh)$.

(iii) For every $x\in M(g)\cap M(h)$, there exists a local inverse pair between $g$ and $h$ containing $x$.

(iv) $|M(gh)|$ is equal to the sum of the numbers of local inverse pairs and free points between $g$ and $h$.

\end{lem}
\demo (i) Let $x\in M(g)$ and $x\in Fix(gh)$. Then we see that there exist two cycle factors $g_i\in\{g\}$ and $h_j\in\{h\}$ such that they can be expressed as $g_i=(x,y,...)$ and $h_j=(y,x,...)$. Since $g_ih_j\neq(1)$, there exists a local inverse pair containing $x$ between $g$ and $h$ by Definition\ \ref{pan2-0}.

(ii) Without loss of generality, we set $\{A,A^{-1}\}$ is a local inverse pair between $g$ and $h$ where $A=x_1,x_2,...,x_l$ and $A^{-1}=x_l,x_{l-1},...,x_1$. According to Definition\ \ref{pan2-0}, if $l>1$, then $x_i\in Fix(gh)$ for $i=1,...,l-1$ and $x_l\notin Fix(gh)$; if $l=1$, $x_l\notin Fix(gh)$.

(iii) Let $x\in M(g)\cap M(h)$. If $x$ is in $Fix(gh)$ or $Fix(hg)$, then by (i) we infer (iii). Otherwise, we have $x^g\neq x^{h^{-1}}$ and $x^{g^{-1}}\neq x^h$. In this case, $\{A,A^{-1}\}$ is a local inverse pair containing $x$ between $g$ and $h$, where $A=x$ and $A^{-1}=x$.

(iv) For all $x\in M(g)\cup M(h)$, $x$ is either in $M(g)\cap M(h)$ or in $(M(g)\cap Fix(h))\cup(Fix(g)\cap M(h))$, and thus $x$ is either a free point or contained in a local inverse pair by Definition\ \ref{pan2-0} and (iii). Then by (i) and (ii) we infer (iv). The proof is completed. \qed

 Using Lemma\ \ref{pan2-2} (i), we obtain the following lemma immediately.

\begin{lem}\label{pan2-1}\normalfont
Let $\alpha=\alpha_1\cdot\cdot\cdot\alpha_s$ and $\beta$ be two permutations in $S_n$. If $x\in M(\alpha^{-1})\cap Fix([\alpha,\beta])$, then there exist $\alpha_i$ and $\alpha_j$ for which one of the following conditions holds:

(i) $\alpha_i^{-1}\alpha_j^\beta=(1)$ and $x\in M(\alpha_i^{-1})$.

(ii) There exists a local inverse pair between $\alpha_i^{-1}$ and $\alpha_j^\beta$ containing $x$.

\end{lem}

So far, we have seen that if there is no transitive cycle factors chain of $\beta$ on $\alpha$, then we can use local inverse pairs and free points to describe $\alpha^{-1}$ and $\alpha^{\beta}$, and so we deal with the case that there is no transitive cycle factors chain of $\beta$ on $\alpha$ in the following.

\begin{lem}\label{pan3-4}\normalfont
Let $\alpha$ and $\beta$ be two permutations in $S_n$ with $1<|M([\alpha,\beta])|\leq4$. If there is no transitive cycle factors chain of $\beta$ on $\alpha$, then $|\{\alpha\}|\leq3$.
\end{lem}
\demo Suppose $\alpha=\alpha_1\alpha_2\cdot\cdot\cdot\alpha_s$ with $\{\alpha\}=\{\alpha_1,\alpha_2,...,\alpha_s\}$ and $s>3$. Note $|M(\alpha_i)|\geq2$ for $i=1,2,...,s$. We have $|M(\alpha)|=|M(\alpha_1)|+|M(\alpha_2)|+\cdot\cdot\cdot+|M(\alpha_s)|\geq8$, and so there exist some points in $M(\alpha)\cup M(\alpha^{\beta})$ that become fixed points of $[\alpha,\beta]$. Since there is no transitive cycle factors chain of $\beta$ on $\alpha$, it follows that $\alpha_i^{-1}\alpha_j^{\beta}\neq(1)$ for all $i$ and $j$. Then by Lemma\ \ref{pan2-2} (i) we deduce that there exist at least one local inverse pair between $\alpha^{-1}$ and $\alpha^{\beta}$, and thus the number of free points between $\alpha^{-1}$ and $\alpha^{\beta}$ is not more than $3$ by Lemma\ \ref{pan2-2} (iv). On the other hand, if there exists a cycle factor $g\in\{\alpha^{-1}\}$ such that there is no local inverse pair between $g$ and $\alpha^{\beta}$, then by Lemma\ \ref{pan2-2} (iii) we see that all elements of $M(g)$ are free points between $\alpha^{-1}$ and $\alpha^{\beta}$, and thus there are at least four free points between $\alpha^{-1}$ and $\alpha^{\beta}$ from the proof of Remark\ \ref{pan3-2046}, a contradiction. Hence, there exist at least one local inverse pair between $g$ and $\alpha^{\beta}$ for each cycle factor $g\in\{\alpha^{-1}\}$, in other words, there exist at least $s$ local inverse pairs between $\alpha^{-1}$ and $\alpha^{\beta}$. However, Lemma\ \ref{pan2-2} (ii) implies that the number of local inverse pairs between $\alpha^{-1}$ and $\alpha^{\beta}$ is not more than $4$, and so $s\leq4$. If $s=4$, then it can only be $\alpha^{-1}=(A)(B)(C)(D)$ and $\alpha^{\beta}=(A^{-1})(B^{-1})(C^{-1})(D^{-1})$, however, $\alpha^{-1}\alpha^{\beta}=(1)$, a contradiction. Thus we have deduced that $s\leq3$.  \qed

\begin{thm}\label{pan3-6}\normalfont
Let $\alpha$ and $\beta$ be two permutations in $S_n$ with $|M([\alpha,\beta])|=4$, where $\alpha$ is a cycle. Then one of the followings holds:

(i) $\alpha^{-1}=(x,y)$ and $\alpha^{\beta}=(x',y')$ where $x,y,x',y'$ are four free points between $\alpha^{-1}$ and $\alpha^{\beta}$.

(ii) $\alpha^{-1}=(ABx)$ and $\alpha^{\beta}=(A^{-1}B^{-1}y)$ where $x,y$ are two free points, $\{A,A^{-1}\}$ and $\{B,B^{-1}\}$ are two local inverse pairs  between $\alpha^{-1}$ and $\alpha^{\beta}$.

(iii) $\alpha^{-1}=(ABCD)$ and $\alpha^{\beta}=(A^{-1}B^{-1}C^{-1}D^{-1})$ where $\{A,A^{-1}\}$, $\{B,B^{-1}\}$, $\{C,C^{-1}\}$, $\{D,D^{-1}\}$ are four local inverse pairs between $\alpha^{-1}$ and $\alpha^{\beta}$.
\end{thm}
\demo It follows from Lemma\ \ref{pan2-2} (iv) that the number of free points is at most $4$. In addition, Remark\ \ref{pan3-2046} shows that the number of free points may be $0$ or $2$ or $4$. Thus we divide into three cases to discuss $\alpha^{-1}$ and $\alpha^{\beta}$.

(i) If there are four free points, then there is no local inverse pair by Lemma\ \ref{pan2-2} (iv), and thus we obtain (i) immediately.

(ii) If there are two free points, then by Lemma\ \ref{pan2-2} (iv) we see that there are two local inverse pairs. So we may assume that the block partitions of $\alpha^{-1}$ and $\alpha^{\beta}$ are $\{A,B,x\}$ and $\{A^{-1},B^{-1},y\}$ respectively, where $x,y$ are two free points, $\{A,A^{-1}\}$ and $\{B,B^{-1}\}$ are two local inverse pairs  between $\alpha^{-1}$ and $\alpha^{\beta}$. Without loss of generality, we take $\alpha^{-1}=(ABx)$. Note that there are two possible cases of the circle permutation on $A^{-1},B^{-1},y$, those are, $(A^{-1}B^{-1}y)$ or $(B^{-1}A^{-1}y)$. However, if $\alpha^{\beta}=(B^{-1}A^{-1}y)$, then $\{A,A^{-1}\}$ and $\{B,B^{-1}\}$ are not local inverse pairs between $\alpha^{-1}$ and $\alpha^{\beta}$, a contradiction. It follows from Definition\ \ref{pan2-0} that $\alpha^{\beta}=(A^{-1}B^{-1}y)$.

(iii) If there is no free point, then there are four local inverse pairs from Lemma\ \ref{pan2-2} (iv), those are, $\{A,A^{-1}\}$, $\{B,B^{-1}\}$, $\{C,C^{-1}\}$ and $\{D,D^{-1}\}$. Similarly, we set $\alpha^{-1}=(ABCD)$. Obviously, there are six possible cases of the circle permutation on $A^{-1},B^{-1},D^{-1},C^{-1}$. However, nothing but $\alpha^{\beta}=(A^{-1}B^{-1}C^{-1}D^{-1})$, for example, if $\alpha^{\beta}=(A^{-1}B^{-1}D^{-1}C^{-1})$ then $\{C,C^{-1}\}$ and $\{D,D^{-1}\}$ are not local inverse pairs between $\alpha^{-1}$ and $\alpha^{\beta}$, a contradiction.  \qed

\begin{thm}\label{pan3-5}\normalfont
Let $\alpha$ and $\beta$ be two permutations in $S_n$ with $|M([\alpha,\beta])|=3$, where $\alpha$ is a cycle. Then one of the followings holds:

(i) $\alpha^{-1}=(Ax)$ and $\alpha^{\beta}=(A^{-1}y)$ where $x,y$ are two free points, $\{A,A^{-1}\}$ is a local inverse pair between $\alpha^{-1}$ and $\alpha^{\beta}$.

(ii) $\alpha^{-1}=(ABC)$ and $\alpha^{\beta}=(A^{-1}B^{-1}C^{-1})$ where $\{A,A^{-1}\}$, $\{B,B^{-1}\}$ and $\{C,C^{-1}\}$ are three local inverse pairs between $\alpha^{-1}$ and $\alpha^{\beta}$.
\end{thm}
\demo It follows from $|M([\alpha,\beta])|=3$ and Remark\ \ref{pan3-2046} that the number of free points may be $0$ or $2$. So we divide into two cases to discuss $\alpha^{-1}$ and $\alpha^{\beta}$.

(i) If there are two free points, then there exists only one local inverse pair $\{A,A^{-1}\}$, and thus we may set $\alpha^{-1}=(Ax)$ and $\alpha^{\beta}=(A^{-1}y)$ where $x,y$ are two free points.

(ii) If there is no free point, then there are three local inverse pairs. So we may set $\alpha^{-1}=(ABC)$. Proceeding as in the proof of Theorem\ \ref{pan3-6} (ii), we derive $\alpha^{\beta}=(A^{-1}B^{-1}C^{-1})$.  \qed

\begin{thm}\label{pan3-7}\normalfont
Let $\alpha$ and $\beta$ be two permutations in $S_n$ with $|\{\alpha\}|=2$, $|M([\alpha,\beta])|=4$ and there is no transitive cycle factors chain of $\beta$ on $\alpha$. Then one of the followings holds:

(i) $\alpha^{-1}=(A)(Bx)$ and $\alpha^{\beta}=(B^{-1})(A^{-1}y)$ where $x,y$ are two free points, $\{A,A^{-1}\}$ and $\{B,B^{-1}\}$ are two local inverse pairs between $\alpha^{-1}$ and $\alpha^{\beta}$.

(ii) Either $\alpha^{-1}=(ABC)(D)$, $\alpha^{\beta}=(D^{-1}B^{-1}C^{-1})(A^{-1})$ or $\alpha^{\beta}=(A^{-1}C^{-1})(B^{-1}D^{-1})$, $\alpha^{-1}=(AB)(CD)$ where $\{A,A^{-1}\}$, $\{B,B^{-1}\}$, $\{C,C^{-1}\}$, $\{D,D^{-1}\}$ are four local inverse pairs between $\alpha^{-1}$ and $\alpha^{\beta}$.

\end{thm}
\demo Since there is no transitive cycle factors chain of $\beta$ on $\alpha$, it follows from Lemma\ \ref{pan2-2} (iv) that the number of free points and local inverse pairs is $4$. Additionally, if there is no local inverse pair, then the number of free points between $\alpha^{-1}$ and $\alpha^{\beta}$ is not less than $8$ because of $|\{\alpha\}|=2$ and the proof of Remark\ \ref{pan3-2046}. Hence, there exists at least one local inverse pair, which means that the number of free points may be $0$ or $2$.

(i) If there are two free points $x$ and $y$, then there are two local inverse pairs $\{A,A^{-1}\}$ and $\{B,B^{-1}\}$, and so we can let $\alpha^{-1}=(A)(Bx)$. Note that for $\alpha^{\beta}$, there are two possible cases which are $(A^{-1})(B^{-1}y)$ and $(B^{-1})(A^{-1}y)$. However, if $\alpha^{\beta}=(A^{-1})(B^{-1}y)$, then $\{A,A^{-1}\}$ is not a local inverse pair between $\alpha^{-1}$ and $\alpha^{\beta}$, a contradiction. Therefore, $\alpha^{\beta}=(B^{-1})(A^{-1}y)$.

(ii) If there is no free point, then there are four local inverse pairs $\{A,A^{-1}\}$, $\{B,B^{-1}\}$, $\{C,C^{-1}\}$ and $\{D,D^{-1}\}$. We observe that there are two possible cases:

1. There exists one cycle factor of $\alpha^{-1}$ has three blocks and the other cycle factor has one block. In this case, an argument similar to the one used in (i) shows that $\alpha^{-1}=(ABC)(D)$ and $\alpha^{\beta}=(D^{-1}B^{-1}C^{-1})(A^{-1})$

2. Every cycle factor of $\alpha^{-1}$ has two blocks. In the same way, we derive $\alpha^{-1}=(AB)(CD)$ and $\alpha^{\beta}=(A^{-1}C^{-1})(B^{-1}D^{-1})$. \qed

\begin{thm}\label{pan3-8}\normalfont
Let $\alpha$ and $\beta$ be two permutations in $S_n$ with $|\{\alpha\}|=2$, $|M([\alpha,\beta])|=3$ and there is no transitive cycle factors chain of $\beta$ on $\alpha$. Then $\alpha^{-1}=(AB)(C)$ and $\alpha^{\beta}=(B^{-1}C^{-1})(A^{-1})$ where $\{A,A^{-1}\}$, $\{B,B^{-1}\}$ and $\{C,C^{-1}\}$ are three local inverse pairs between $\alpha^{-1}$ and $\alpha^{\beta}$.
\end{thm}
\demo According to the proof of Theorem\ \ref{pan3-7}, it follows that the number of free points is $0$ or $2$. If the number of free points is $2$, then there exists only one local inverse pair, and thus there exists a cycle factor $\alpha_i\in\{\alpha^{-1}\}$ such that there is no local inverse pair between $\alpha_i^{-1}$ and $\alpha^{\beta}$, which implies that the number of free points between $\alpha^{-1}$ and $\alpha^{\beta}$ is not less than $4$ because of $|M(\alpha_i^{-1})|\geq2$ and the proof of Remark\ \ref{pan3-2046}, a contradiction. Hence, there is no free point between $\alpha^{-1}$ and $\alpha^{\beta}$, and then by using the way of proving Theorem\ \ref{pan3-7} to derive this theorem. \qed

\begin{thm}\label{pan3-20}\normalfont
Let $\alpha$ and $\beta$ be two permutations in $S_n$ with $|\{\alpha\}|=3$, $1<|M([\alpha,\beta])|\leq4$ and there is no transitive cycle factors chain of $\beta$ on $\alpha$. Then $\alpha^{-1}=(AB)(C)(D)$ and $\alpha^{\beta}=(C^{-1}D^{-1})(A^{-1})(B^{-1})$ where $\{A,A^{-1}\}$, $\{B,B^{-1}\}$, $\{C,C^{-1}\}$, $\{D,D^{-1}\}$ are four local inverse pairs between $\alpha^{-1}$ and $\alpha^{\beta}$.
\end{thm}
\demo An argument similar to the one used in Theorem\ \ref{pan3-8} shows that there is no free points between $\alpha^{-1}$ and $\alpha^{\beta}$, and there are four local inverse pairs $\{A,A^{-1}\}$, $\{B,B^{-1}\}$, $\{C,C^{-1}\}$ and $\{D,D^{-1}\}$ between $\alpha^{-1}$ and $\alpha^{\beta}$. Proceeding as in the proof of Theorem\ \ref{pan3-7}, we prove this theorem. \qed

We conclude this section by pointing out two useful remarks.

\begin{Remark}\label{pan3-206}\normalfont
In Theorems\ \ref{pan3-6}~-\ \ref{pan3-20}, we see that there exists at most one point that is in $M(\alpha^\beta)$ but not in $M(\alpha)$ except Theorem\ \ref{pan3-6} (i).
\end{Remark}

\begin{Remark}\label{pan3-1206}\normalfont
In Theorems\ \ref{pan3-7}~-\ \ref{pan3-20}, we see that there are at most two cycle factors of $\alpha$ whose lengths are equal, those are, Theorem\ \ref{pan3-7} (ii) 2 and Theorem\ \ref{pan3-20}.
\end{Remark}

\section {Simultaneous conjugation}

In this section, we first prove Conjecture\ \ref{pan1-0} in the case when there is no transitive cycle factors chain of $\beta$ on $\alpha$. Note that for the case that there is no transitive cycle factors chain of $\beta$ on $\alpha$, it suffices to prove each case of Theorems\ \ref{pan3-6}~-\ \ref{pan3-20}. However, there are still two problems that need to be solved, one is that for a fixed $\alpha$, there exist many $\beta'$ such that $[\alpha,\beta']=[\alpha,\beta]$, and the other is that we can't give an exact form of $\beta$ even if $\alpha$ and $\alpha^\beta$ have been given. In order to settle these two problems, we give the following two lemmas.

\begin{lem}\label{pan2-66}\normalfont
Let $\alpha=\alpha_1\cdot\cdot\cdot\alpha_k$ and $\beta$ be two permutations in $S_n$ with $|M(\beta)\setminus M(\alpha)|\leq1$ and there are at most two cycle factors of $\alpha$ whose lengths are equal. If there exists an involution $\omega$ such that $M(\omega)\subseteq M(\alpha)$ and $\alpha_i^\omega=\alpha_i^{-1}$ for $i=1,2,...,k$ and $\beta^\omega=\beta^{-1}$, then the followings hold:

(i) If all lengths of the cycle factors of $\alpha$ are different from each other, then for each $\beta'$ with $\alpha^{\beta'}=\alpha^{\beta}$, there exists a permutation $\gamma$ such that $\beta'^\gamma=\beta'^{-1}$ and $\alpha_i^\gamma=\alpha_i^{-1}$ for $i=1,2,...,k$.

(ii) If $|\alpha_1|=|\alpha_2|$ and there exists an involution $\mu\in C_{S_n}(\alpha)\cap C_{S_n}(\beta)\cap C_{S_n}(\omega)$ such that $\alpha_1^\mu=\alpha_2$ and $\alpha_2^\mu=\alpha_1$ and $M(\mu)\subseteq M(\alpha)$, then for each $\beta'$ with $\alpha^{\beta'}=\alpha^{\beta}$, there exists a permutation $\gamma$ such that $\beta'^\gamma=\beta'^{-1}$ and $\alpha_i^\gamma=\alpha_i^{-1}$ for $i=1,2,...,k$.
\end{lem}
\demo (i) It is well known that $\beta'\in C_{S_n}(\alpha)\beta$ for each $\beta'$ with $\alpha^{\beta'}=\alpha^{\beta}$. Since all lengths of the cycle factors of $\alpha$ are different from each other, it follows that $$C_{S_n}(\alpha)=\langle\alpha_1\rangle\times\langle\alpha_2\rangle\cdot\cdot\cdot\times\langle\alpha_k\rangle\times S_{([n]\setminus M(\alpha))},{\rm{where}}~ S_{([n]\setminus M(\alpha))}~{\rm{is~ the~ symmetric~ group~ on}}~[n]\setminus M(\alpha).$$
So we may set $\beta'=ab\beta$ where $a=\alpha_1^{i_1}\alpha_2^{i_2}\cdot\cdot\cdot\alpha_k^{i_k}$ and $b\in S_{([n]\setminus M(\alpha))}$. As with $|M(\beta)\setminus M(\alpha)|\leq1$ and Lemma\ \ref{pan2-7} (i), there exists an involution $\theta\in S_{([n]\setminus M(\alpha))}$ such that $\beta^\theta=\beta$ and $b^\theta=b^{-1}$. Picking $\gamma=ab\omega\theta$. Note $\beta'=\beta^{(ab)^{-1}}ab$ and $ab=ba$. We see $\alpha_i^\gamma=\alpha_i^{\omega\theta}=(\alpha_i^{-1})^\theta=\alpha_i^{-1}$ for $i=1,2,...,k$, and  $\beta'^{\gamma}=(\beta^{(ab)^{-1}}ab)^{ab\omega\theta}=\beta^{\omega\theta}(ab)^{\omega\theta}=\beta^{-1}a^{-1}b^{-1}=\beta^{-1}(ab)^{-1}=\beta'^{-1}$, as desired.

(ii) Note that $C_{S_n}(\alpha)=\langle\alpha_3\rangle\times\cdot\cdot\cdot\times\langle\alpha_k\rangle\times G\times S_{([n]\setminus M(\alpha))}$ where $G=(\langle\alpha_1\rangle\times\langle\alpha_2\rangle)\rtimes\langle\mu\rangle$. So we may set $\beta'=ab\beta$ where $a=\alpha_1^{i_1}\alpha_2^{i_2}\cdot\cdot\cdot\alpha_k^{i_k}\mu^i$ and $b\in S_{([n]\setminus M(\alpha))}$. If $i=0$, the (ii) follows by (i). Consider $i=1$. Then by $|M(\beta)\setminus M(\alpha)|\leq1$ and Lemma\ \ref{pan2-7} (i), there exists an involution $\theta\in S_{([n]\setminus M(\alpha))}$ such that $b^\theta=b^{-1}$ and $\beta^\theta=\beta$. Choose
$\gamma=ab\omega\mu\theta$. Note $\beta'=\beta^{(ab)^{-1}}ab$ and $ab=ba$. We see $\alpha_i^\gamma=\alpha_i^{-1}$ for $i=1,2,...,k$, and $\beta'^{\gamma}=(\beta^{(ab)^{-1}}ab)^{ab\omega\mu\theta}=\beta^{\omega\mu\theta}(ab)^{\omega\mu\theta}=\beta^{-1}(\alpha_1^{-i_1}\alpha_2^{-i_2}\cdot\cdot\cdot\alpha_k^{-i_k}\mu b)^{\mu\theta}=\beta^{-1}[\mu(\alpha_1^{-i_1}\alpha_2^{-i_2}\cdot\cdot\cdot\alpha_k^{-i_k}b)^\mu]^{\mu\theta}=\beta^{-1}(ab)^{-1}=\beta'^{-1}$, as desired. \qed

\begin{lem}\label{pan2-12}\normalfont
Let $\alpha$ and $\beta$ be two permutations in $S_n$ with $M(\beta)\subseteq M(\alpha)\cup M(\alpha^\beta)$ and there exists at most one point that is in $M(\alpha^\beta)$ but not in $M(\alpha)$. If there exists an involution $\omega$ such that $M(\omega)\subseteq M(\alpha)$, $\alpha^\omega=\alpha^{-1}$ and $x^{\omega\beta\omega\beta}=x$ for all $x\in M(\alpha)$, then $\beta^\omega=\beta^{-1}$.
\end{lem}
\demo Suppose there exists at most one point which is in $M(\alpha^\beta)$ but not in $M(\alpha)$. Then by Assumption 2.5 and $M(\beta)\subseteq M(\alpha)\cup M(\alpha^\beta)$ we see that there exists at most one point which is in $M(\beta)$ but not in $M(\alpha)$, and thus the $x^{\omega\beta\omega\beta}=x$ for all $x\in M(\alpha)$ implies  $\omega\beta\omega\beta$ is the identity, and therefore $\omega\beta\omega=\beta^{-1}$. Additionally, since $\omega$ is an involution, we derive $\omega\beta\omega=\beta^\omega=\beta^{-1}$. The proof of this lemma is now complete.   \qed

According to Remark\ \ref{pan3-206}, Remark\ \ref{pan3-1206} and Lemma\ \ref{pan2-66}, it suffices to consider one $\beta$ with $M(\beta)\subseteq M(\alpha)\cup M(\alpha^\beta)$ for each case of Theorems\ \ref{pan3-6}~-\ \ref{pan3-20} except Theorem\ \ref{pan3-6} (i). Now we start to prove Conjecture\ \ref{pan1-0} in the case when there is no transitive cycle factors chain of $\beta$ on $\alpha$.

\begin{Pro}\label{pan4-1116}\normalfont
Let $\alpha$ and $\beta$ be two permutations in $S_n$ with $0<|M([\alpha,\beta])|\leq4$, where $\alpha$ is a cycle. Then there exists a permutation $\gamma$ such that $\alpha^{\gamma}=\alpha^{-1}$ and $\beta^{\gamma}=\beta^{-1}$.
\end{Pro}
\demo Apparently, it suffices to discuss each case of Theorem\ \ref{pan3-6} and Theorem\ \ref{pan3-5}, as follows:

Case 1: Theorem\ \ref{pan3-6} (i). Without loss of generality, we set $x^\beta=x'$ and $y^\beta=y'$. Note that $x$ and $x'$ are contained in the same cycle factor of $\beta$ as well as $y$ and $y'$. By Lemma\ \ref{pan2-7} (i) (ii), we see that there exists an involution $\gamma$ such that $\alpha^\gamma=\alpha^{-1}$ and $\beta^\gamma=\beta^{-1}$.

Case 2: Theorem\ \ref{pan3-6} (ii). According to Definition\ \ref{pan2-0} and Lemma\ \ref{pan2-66}, we may assume that
$$\alpha^{-1}=(x_1,x_2,...,x_r)~{\rm{and}}~\alpha^{\beta}=(x_{s},x_{s-1},...,x_1,x_{r-1},...,x_{s+1},x_{r+1}),$$
where $1\leq s\leq r-2$ and $x_s^\beta=x_{s},x_{s-1}^\beta=x_{s-1},...,x_r^\beta=x_{r-1},x_{r-1}^\beta=x_{r-2},...,x_{s+2}^\beta=x_{s+1},x_{s+1}^\beta=x_{r+1}$ and $M(\beta)\subseteq M(\alpha)\cup M(\alpha^\beta)$. Then we have $\beta=(x_r,x_{r-1},...,x_{s+1},x_{r+1})$. Pick an involution $\omega$ such that $x_1^\omega=x_s,~x_2^\omega=x_{s-1},...,x_s^\omega=x_{1},~x_{s+1}^\omega=x_{r},~x_{s+2}^\omega=x_{r-1},...,x_r^\omega=x_{s+1}$ and $M(\omega)\subseteq M(\alpha)$. One easily checks that $\alpha^\omega=\alpha^{-1}$ and $\beta^\omega=\beta^{-1}$, as desired.

Case 3: Theorem\ \ref{pan3-6} (iii). Similarly, we may assume that $\alpha^{-1}=(x_1,x_2,...,x_r)$ and $$\alpha^{\beta}=(x_{s},x_{s-1},...,x_1,x_t,x_{t-1},...,x_{s+1},x_k,x_{k-1},...,x_{t+1},x_r,x_{r-1},...,x_{k+1}),$$
where $1\leq s<t<k<r$ and $x_r^\beta=x_{s},x_{r-1}^\beta=x_{s-1},...,x_{r-s+1}^\beta=x_{1},x_{r-s}^\beta=x_{t},...,x_{r-t+1}^\beta=x_{s+1},x_{r-t}^\beta=x_{k},...,x_{r-k+1}^\beta=x_{t+1},x_{r-k}^\beta=x_{r},...,x_{1}^\beta=x_{k+1}$ and $M(\beta)\subseteq M(\alpha)$. Note that
\[
\begin{cases}
  x_{j}^\beta=x_{j+s-r}, \ &\text{if $r-s+1\leq j\leq r$ \ }\\
  x_{j}^\beta=x_{s+j+t-r}, \ &\text{if $r-t+1\leq j\leq r-s$ \ }\\
 x_{j}^\beta=x_{j+t+k-r}, \ &\text{if $r-k+1\leq j\leq r-t$ \ }\\
    x_{j}^\beta=x_{j+k},\ &\text{if $1\leq j\leq r-k$}
\end{cases}
\]
Take an involution $\omega$ so that $x_1^\omega=x_r,x_2^\omega=x_{r-1},...,x_i^\omega=x_{r-i+1},...,x_r^\omega=x_{1}$ and $M(\omega)\subseteq M(\alpha)$. Obviously, $\alpha^\omega=\alpha^{-1}$. On the other hand, an easy computation to verify that
\[
 x_j^{\omega\beta\omega\beta}=x_{r-j+1}^{\beta\omega\beta}=
\begin{cases}
 x_{s-j+1}^{\omega\beta} = x_{r+j-s}^{\beta}=x_j, \ &\text{if $1\leq j\leq s$ \ }\\
   x_{s+t-j+1}^{\omega\beta}= x_{r+j-s-t}^{\beta}=x_j, \ &\text{if $s+1\leq j\leq t$}\\
x_{k+t-j+1}^{\omega\beta}= x_{r+j-k-t}^{\beta}=x_j, \ &\text{if $t+1\leq j\leq k$}\\
x_{k+r-j+1}^{\omega\beta}= x_{j-k}^{\beta}=x_j, \ &\text{if $k+1\leq j\leq r$}
\end{cases}
\]
As with Lemma\ \ref{pan2-12}, as desired.

Case 4: Theorem\ \ref{pan3-5} (i). In the same way, we take $\alpha^{-1}=(x_1,x_2,...,x_{r-1},x)$ and $\alpha^{\beta}=(x_{r-1},...,x_1,y)$, where $x_{r-1}^\beta=x_{r-1},x_{r-2}^\beta=x_{r-2},...,x_1^\beta=x_{1},x^\beta=y$ and $M(\beta)\subseteq M(\alpha)\cup M(\alpha^\beta)$. Then we have $\beta=(x,y)$.
Regard an involution $\omega$ with $x_1^\omega=x_{r-1},~x_2^\omega=x_{r-2},...,x_1^\omega=x_{r-1},~x^\omega=x$ and $M(\omega)\subseteq M(\alpha)$, as desired.

Case 4: Theorem\ \ref{pan3-5} (ii). Similarly, we may set that $$\alpha^{-1}=(x_1,x_2,...,x_r)~{\rm{and}}~\alpha^{\beta}=(x_{s},x_{s-1},...,x_1,x_t,x_{t-1},...,x_{s+1},x_r,x_{r-1},...,x_{t+1}),$$
where $1\leq s<t<r$ and $x_t^\beta=x_{s},x_{t-1}^\beta=x_{s-1},...,x_{t-s+1}^\beta=x_{1},x_{t-s}^\beta=x_{t},...,x_{1}^\beta=x_{s+1},x_{r}^\beta=x_{r},...,x_{t+1}^\beta=x_{t+1}$ and $M(\beta)\subseteq M(\alpha)$. Note that
\[
\begin{cases}
  x_{j}^\beta=x_{j}, \ &\text{if $t+1\leq j\leq r$ \ }\\
  x_{j}^\beta=x_{s+j-t}, \ &\text{if $t-s+1\leq j\leq t$ \ }\\
 x_{j}^\beta=x_{j+s}, \ &\text{if $1\leq j\leq t-s$ \ }
\end{cases}
\]
Take an involution $\omega$ so that $x_1^\omega=x_t,x_2^\omega=x_{t-1},...,x_{t+1}^\omega=x_{r},...,x_r^\omega=x_{t+1}$ and $M(\omega)\subseteq M(\alpha)$. Obviously, $\alpha^\omega=\alpha^{-1}$. On the other hand, one easily checks that
\[
 x_j^{\omega\beta\omega\beta}=
\begin{cases}
 x_{t-j+1}^{\beta\omega\beta}=x_{s-j+1}^{\omega\beta} = x_{t+j-s}^{\beta}=x_j, \ &\text{if $1\leq j\leq s$ \ }\\
   x_{t-j+1}^{\beta\omega\beta}=x_{s+t-j+1}^{\omega\beta}= x_{j-s}^{\beta}=x_j, \ &\text{if $s+1\leq j\leq t$}\\
x_{r-j+t+1}^{\beta\omega\beta}=x_{r+t-j+1}^{\omega\beta}= x_{j}^{\beta}=x_j, \ &\text{if $t+1\leq j\leq r$}
\end{cases}
\]
As with Lemma\ \ref{pan2-12}, as desired.    \qed

\begin{Pro}\label{pan4-36}\normalfont
Let $\alpha=\alpha_1\alpha_2$ and $\beta$ be in $S_n$ with $\alpha_i^\beta\notin\{\alpha\}$ for $i=1,2$, where $\alpha_1$ and $\alpha_2$ are two cycle factors of $\alpha$. If $0<|M([\alpha,\beta])|\leq4$, then there exists a permutation $\gamma\in S_n$ such that $\alpha_1^{\gamma}=\alpha_1^{-1}$ and $\alpha_2^{\gamma}=\alpha_2^{-1}$ and $\beta^{\gamma}=\beta^{-1}$.
\end{Pro}
\demo Note that it is suffices to prove each case of Theorem\ \ref{pan3-7} and Theorem \ref{pan3-8}, as follows:

Case 1: Theorem\ \ref{pan3-7} (i). According to Definition\ \ref{pan2-0} and Lemma\ \ref{pan2-66}, we can assume that $$\alpha^{-1}=(x_1,...,x_r)(y_1,...,y_{l})~{\rm{and}}~\alpha^{\beta}=(y_{l},...,y_{1},x_{r+1})(x_{r-1},x_{r-2},...,x_1),$$ where $x_l^\beta=y_{l},x_{l-1}^\beta=x_{l-1},...,x_{1}^\beta=y_{1},~x_{r}^\beta=x_{r+1}$ and $y_l^\beta=x_{l},y_{l-1}^\beta=x_{l-1},..,y_{1}^\beta=x_{1}$ and $M(\beta)\subseteq M(\alpha)\cup M(\alpha^\beta)$. Then we see $\beta=(x_1,y_1)(x_2,y_2)\cdot\cdot\cdot(x_l,y_l)(x_r,x_{r+1})$. Choosing an involution $\omega$ such that $x_{l}^{\omega}=x_1,x_{l-1}^{\omega}=x_2,...,x_1^{\omega}=x_{l}, x_r^{\omega}=x_{r}$ and $y_l^{\omega}=y_1,y_{l-1}^{\omega}=y_2,...,y_1^{\omega}=y_{l}$ and $M(\omega)\subseteq M(\alpha)$. One easily checks that $\alpha_1^\omega=\alpha_1^{-1}$, $\alpha_2^\omega=\alpha_2^{-1}$ and $\beta^\omega=\beta^{-1}$, as desired.

Case 2: $\alpha^{-1}=(ABC)(D)$ and $\alpha^{\beta}=(D^{-1}B^{-1}C^{-1})(A^{-1})$ of Theorem\ \ref{pan3-7} (ii). We may let  $$\alpha^{-1}=(x_1,...,x_r)(y_1,...,y_{l})~{\rm{and}}~\alpha^{\beta}=(y_{l},...,y_{1},x_s,x_{s-1},...,x_{l+1},x_r,...,x_{s+1})(x_{l},x_{l-1},...,x_1),$$
where $1<l<s<r$ and $x_r^\beta=y_l,x_{r-1}^\beta=y_{l-1},...,~x_{1}^\beta=x_{s+1}$ and $y_l^\beta=x_l,y_{l-1}^\beta=x_{l-1},...,y_1^\beta=x_1$ and $M(\beta)\subseteq M(\alpha)$. It is easy to see that
\[
\begin{cases}
  x_{j}^\beta=y_{j+l-r}, \ &\text{if $r-l+1\leq j\leq r$ \ }\\
  x_{j}^\beta=x_{s+j+l-r}, \ &\text{if $r-s+1\leq j\leq r-l$ \ }\\
    x_{j}^\beta=x_{j+s},\ &\text{if $1\leq j\leq r-s$}
\end{cases}
\]
Pick an involution $\omega$ such that $x_1^{\omega}=x_r,x_2^{\omega}=x_{r-1},...,x_r^{\omega}=x_{1}$ and $y_1^{\omega}=y_l,y_2^{\omega}=y_{l-1},...,y_l^{\omega}=y_{1}$ and $M(\omega)\subseteq M(\alpha)$. Note $y_i^{\omega\beta\omega\beta}=y_{l-i+1}^{\beta\omega\beta}=x_{l-i+1}^{\omega\beta}=x_{r-l+i}^\beta=y_i$ for $i=1,...,l$ and further
\[
 x_j^{\omega\beta\omega\beta}=x_{r-j+1}^{\beta\omega\beta}=
\begin{cases}
 y_{l-j+1}^{\omega\beta} = y_{j}^{\beta}=x_j, \ &\text{if $1\leq j\leq l$ \ }\\
   x_{s+l-j+1}^{\omega\beta}= x_{r+j-s-l}^{\beta}=x_j, \ &\text{if $l+1\leq j\leq s$}\\
x_{r+s-j+1}^{\omega\beta}= x_{j-s}^{\beta}=x_j, \ &\text{if $s+1\leq j\leq r$}
\end{cases}
\]
As with Lemma\ \ref{pan2-12}, as desired.

Case 3: $\alpha^{-1}=(AB)(CD)$, $\alpha^{\beta}=(A^{-1}C^{-1})(B^{-1}D^{-1})$ of Theorem\ \ref{pan3-7} (ii). We may set $$\alpha^{-1}=(x_1,...,x_r)(y_1,...,y_{l})~{\rm{and}}~\alpha^{\beta}=(x_{r},x_{r-1},...,x_{s+1},y_{s},...,y_{1})(y_l,...,y_{s+1},x_s,...,x_1),$$
where $1\leq s<l\leq r$ and $x_r^\beta=x_{r},x_{r-1}^\beta=x_{r-1},...,x_{s+1}^\beta=x_{s+1},~x_{s}^\beta=y_s,...,x_{1}^\beta=y_1$ and $y_l^\beta=y_{l},y_{l-1}^\beta=y_{l-1},..,y_{s+1}^\beta=y_{s+1},y_{s}^\beta=x_{s},...,y_{1}^\beta=x_{1}$ and $M(\beta)\subseteq M(\alpha)$. We see that $\beta=(x_s,y_s)\cdot\cdot\cdot(x_1,y_1)$. Take an involution $\omega$ such that $x_1^\omega=x_s,x_2^\omega=x_{s-1},...,x_s^\omega=x_1,x_{s+1}^\omega=x_{r},...,x_{r}^\omega=x_{s+1}$ and $y_1^\omega=y_s,y_2^\omega=y_{s-1},...,y_s^\omega=y_1,y_{s+1}^\omega=y_{l},...,y_{l}^\omega=y_{s+1}$ and $M(\omega)\subseteq M(\alpha)$. It is easy to check that $\alpha_1^{\omega}=\alpha_1^{-1}$ and $\alpha_2^{\omega}=\alpha_2^{-1}$ and $\beta^{\omega}=\beta^{-1}$. In addition, we observe that if $r=l$, then the involution $\mu=(x_r,y_r)(x_{r-1},y_{r-1})\cdot\cdot\cdot(x_1,y_1)$ such that $\mu\in C_{S_n}(\alpha)\cap C_{S_n}(\beta)\cap C_{S_n}(\omega)$ and $\alpha_1^\mu=\alpha_2$ and $\alpha_2^\mu=\alpha_1$, as desired.

Case 4: Theorem \ref{pan3-8}. In this case, we may assume that $$\alpha^{-1}=(x_1,...,x_r)(y_1,...,y_{l})~{\rm{and}}~\alpha^{\beta}=(y_l,y_{l-1},...,y_{1},x_r,...,x_{l+1})(x_{l},x_{l-1},...,x_1),$$
where $1<l<r$ and $x_l^\beta=y_{l},x_{l-1}^\beta=y_{l-1},...,~x_{1}^\beta=x_{1},x_r^\beta=x_r,...,x_{l+1}^\beta=x_{l+1}$ and $y_l^\beta=x_{l},y_{l-1}^\beta=x_{l-1},..,y_{1}^\beta=x_{1}$ and $M(\beta)\subseteq M(\alpha)$. Obviously, $\beta=(x_1,y_1)(x_2,y_2)\cdot\cdot\cdot(x_l,y_l)$.
Pick an involution $\omega$ such that $x_1^{\omega}=x_l,x_2^{\omega}=x_{l-1},...,x_l^{\omega}=x_1,x_{l+1}^{\omega}=x_r,...,x_r^{\omega}=x_{l+1}$ and
$y_1^{\omega}=y_{l},y_2^{\omega}=y_{l-1},...,y_l^{\omega}=y_{1}$ and $M(\omega)\subseteq M(\alpha)$. It is straightforward to show $\alpha_1^{\omega}=\alpha_1^{-1}$ and $\alpha_2^{\omega}=\alpha_2^{-1}$ and $\beta^{\omega}=\beta^{-1}$.   \qed

\begin{Pro}\label{pan4-11}\normalfont
Let $\alpha=\alpha_1\alpha_2\alpha_3$ and $\beta$ be in $S_n$ with $\alpha_i^\beta\notin\{\alpha\}$ for $i=1,2,3$, where $\alpha_1,\alpha_2,\alpha_3$ are three cycle factors of $\alpha$. If $0<|M([\alpha,\beta])|\leq4$, then there exists a permutation $\gamma\in S_n$ such that $\alpha_1^{\gamma}=\alpha_1^{-1}$ and $\alpha_2^{\gamma}=\alpha_2^{-1}$ and $\alpha_3^{\gamma}=\alpha_3^{-1}$ and $\beta^{\gamma}=\beta^{-1}$.
\end{Pro}
\demo As with Theorem\ \ref{pan3-20} and Definition\ \ref{pan2-0} and Lemma\ \ref{pan2-66}, we set $\alpha_1=(x_1,x_2,...,x_r)$, $\alpha_2=(y_1,y_2,...,y_s)$, $\alpha_3=(z_1,...,z_t)$ and $$\alpha^{\beta}=(y_s,y_{s-1}...,y_1,z_t,z_{t-1},...,z_1)(x_s,x_{s-1},...,x_1)(x_r,x_{r-1},...,x_{s+1}),$$
where $x_s^\beta=y_s,x_{s-1}^\beta=y_{s-1},...,x_{1}^\beta=y_1,x_{r}^\beta=z_t,...,x_{s+1}^\beta=z_1$ and $y_s^\beta=x_s,y_{s-1}^\beta=x_{s-1},...,y_{1}^\beta=x_1$ and $z_t^\beta=x_r,z_{t-1}^\beta=x_{r-1},...,z_1^\beta=x_{s+1}$ and $M(\beta)\subseteq M(\alpha)$. In this case, we note that $$\beta=(x_1,y_1)(x_2,y_2)\cdot\cdot\cdot(x_s,y_s)(x_{s+1},z_1)\cdot\cdot\cdot(x_r,z_t).$$ Take an involution $\omega$ such that $x_1^{\omega}=x_s,x_2^{\omega}=x_{s-1},...,x_r^{\omega}=x_{s+1}$ and $y_1^{\omega}=y_s,y_{2}^{\omega}=y_{s-1},...,y_s^{\omega}=y_{1}$ and $z_1^{\omega}=z_t,z_2^{\omega}=x_{t-1},...,z_t^{\omega}=z_{1}$ and $M(\omega)\subseteq M(\alpha)$. An easy computation to show that $\alpha_1^{\omega}=\alpha_1^{-1}$ and $\alpha_2^{\omega}=\alpha_2^{-1}$ and $\alpha_3^{\omega}=\alpha_3^{-1}$ and $\beta^{\gamma}=\beta^{-1}$.
On the other hand, we observe that if $s=t$, then the involution $\mu=(x_r,x_s)(y_s,z_s)(x_{r-1},x_{s-1})(y_{s-1},z_{s-1})\cdot\cdot\cdot(x_{s+1},x_1)(y_1,z_1)$ such that $\mu\in C_{S_n}(\alpha)\cap C_{S_n}(\beta)\cap C_{S_n}(\omega)$ and $\alpha_2^\mu=\alpha_3$ and $\alpha_3^\mu=\alpha_2$, as desired.  \qed

Using Proposition\ \ref{pan4-1116} and Proposition\ \ref{pan4-36} and Proposition\ \ref{pan4-11}, we obtain the following lemma.

\begin{lem}\label{pan4-00400}\normalfont
Let $\alpha$ and $\beta$ be two permutations in $S_n$ with $0<M([\alpha,\beta])\leq4$ and there is no transitive cycle factors chain of $\beta$ on $\alpha$. Then there exists a permutation $\gamma\in S_n$ such that $\beta^\gamma=\beta^{-1}$ and $\alpha_i^\gamma=\alpha_i^{-1}$ for all $\alpha_i\in\{\alpha\}$.
\end{lem}

Up to now, we have proved the case that there is no transitive cycle factors chain of $\beta$ on $\alpha$. We are now turning to consider the opposite case, and we first give an useful lemma.

\begin{lem}\label{pan4-00}\normalfont
Let $\alpha$ and $\beta$ be two permutations in $S_n$ with $0<M([\alpha,\beta])\leq4$. Then there exist at most three transitive cycle factors chains of $\beta$ on $\alpha$.
\end{lem}
\demo Proof by contradiction. Assume that there exist $s$ transitive cycle factors chains of $\beta$ on $\alpha$ with $s>3$, those are, $\alpha_{m_j+1}\rightarrow\alpha_{m_j+2}\rightarrow\cdot\cdot\cdot\rightarrow\alpha_{m_j+{k_j}}$ for $j=1,2,...,s$. For each transitive cycle factors chain $\alpha_{m_j+1}\rightarrow\alpha_{m_j+2}\rightarrow\cdot\cdot\cdot\rightarrow\alpha_{m_j+{k_j}}$, by using Lemma\ \ref{pan4-004}, we see there exists a permutation $\mu_j\in C_{S_n}(\alpha_{m_j+1}\alpha_{m_j+2}\cdot\cdot\cdot\alpha_{m_j+{k_j}})$ such that $o(\mu_j)=k_j$, $M(\mu_j)=M(\alpha_{m_j+1}\alpha_{m_j+2}\cdot\cdot\cdot\alpha_{m_j+{k_j}})$ and $M(\alpha_{m_j+2}\cdot\cdot\cdot\alpha_{m_j+{k_j}})\subseteq Fix(\mu_j\beta)$. As with Lemma\ \ref{pan4-0004}, we see that $[\alpha,\beta]=[\alpha'\alpha_{m_1+1},\mu_1\beta]$ where $\alpha'\in S_n$ with $\{\alpha'\}=\{\alpha\}\setminus\{\alpha_{m_1+1},\alpha_{m_1+2},...,\alpha_{m_1+{k_1}}\}$. In addition, it follows from $\mu_1\in C_{S_n}(\alpha_{m_1+1}\alpha_{m_1+2}\cdot\cdot\cdot\alpha_{m_1+{k_1}})$ and $M(\mu_1)=M(\alpha_{m_1+1}\alpha_{m_1+2}\cdot\cdot\cdot\alpha_{m_1+{k_1}})$ and Remark\ \ref{pan2-000006} that there exist $s-1$ transitive cycle factors chains of $\mu_1\beta$ on $\alpha'\alpha_{m_1+1}$, those are, $$\alpha_{m_j+1}\rightarrow\alpha_{m_j+2}\rightarrow\cdot\cdot\cdot\rightarrow\alpha_{m_j+{k_j}}~{\rm{for}} ~j=2,...,s.$$ Using Lemma\ \ref{pan4-0004} again, we infer $[\alpha,\beta]=[\alpha'\alpha_{m_1+1},\mu_1\beta]=[\alpha''\alpha_{m_1+1}\alpha_{m_2+1},\mu_1\mu_2\beta]$ where $\alpha''\in S_n$ with $\{\alpha''\}=\{\alpha'\alpha_{m_1+1}\}\setminus\{\alpha_{m_2+1},\alpha_{m_2+2},...,\alpha_{m_2+{k_2}}\}$. Applying Lemma\ \ref{pan4-0004} repeatedly, and finally we get $[\alpha,\beta]=[\alpha'''\alpha_{m_1+1}\alpha_{m_2+1}\cdot\cdot\cdot\alpha_{m_s+1},\mu_1\mu_2\cdot\cdot\cdot\mu_s\beta]$ where $\alpha'''\in S_n$ with $\{\alpha'''\}=\{\alpha\}\setminus(\{\alpha_{m_1+1}\alpha_{m_1+2}\cdot\cdot\cdot\alpha_{m_1+{k_1}}\}\cup\{\alpha_{m_2+1}\alpha_{m_2+2}\cdot\cdot\cdot\alpha_{m_2+{k_2}}\}\cup\cdot\cdot\cdot\cup\{\alpha_{m_s+1}\alpha_{m_s+2}\cdot\cdot\cdot\alpha_{m_s+{k_s}}\})$.
According to Definition\ \ref{pan2-000000}, it follows that there is no transitive cycle factors chain of $\mu_1\mu_2\cdot\cdot\cdot\mu_s\beta$ on $\alpha'''\alpha_{m_1+1}\alpha_{m_2+1}\cdot\cdot\cdot\alpha_{m_s+1}$, however, $s>3$ contradicts Lemma\ \ref{pan3-4}. \qed

In what follows, we will confirm Conjecture\ \ref{pan1-0} in the case when there exist transitive cycle factors chains of $\beta$ on $\alpha$.

\begin{lem}\label{pan4-0}\normalfont
Let $\alpha$ and $\beta$ be two permutations in $S_n$ with $0<M([\alpha,\beta])\leq4$ and there exist transitive cycle factors chains of $\beta$ on $\alpha$. Then there exists a permutation $\gamma\in S_n$ such that $\alpha^\gamma=\alpha^{-1}$ and $\beta^\gamma=\beta^{-1}$.
\end{lem}
\demo Without loss of generality, let $\alpha_1\rightarrow\alpha_2\rightarrow\cdot\cdot\cdot\rightarrow\alpha_k$ be a transitive cycle factors chain of $\beta$ on $\alpha$. Then by Lemma\ \ref{pan4-004} and Lemma\ \ref{pan4-0004}, there exists a permutation $\mu\in C_{S_n}(\alpha_1\alpha_2\cdot\cdot\cdot\alpha_k)$ such that $o(\mu)=k$, $M(\mu)=M(\alpha_1\alpha_2\cdot\cdot\cdot\alpha_k)$, $M(\alpha_2\alpha_3\cdot\cdot\cdot\alpha_k)\subseteq Fix(\mu\beta)$ and $[\alpha,\beta]=[\alpha'\alpha_1,\mu\beta]$, where $\alpha'$ is a permutation with $\{\alpha'\}=\{\alpha\}\setminus\{\alpha_1,\alpha_2,...,\alpha_k\}$. On the other hand, it follows from Lemma\ \ref{pan4-00} that the number of transitive cycle factors chains of $\beta$ on $\alpha$ is $1$ or $2$ or $3$, and so we divide into three cases to confirm this lemma in the following.

Case 1: there exist one transitive cycle factors chain of $\beta$ on $\alpha$. In this case, we note that there is no transitive cycle factors chain of $\mu\beta$ on $\alpha'\alpha_1$, and then by Lemma\ \ref{pan4-00400} and Lemma\ \ref{pan2-8} we deduce that there exists a permutation $\omega$ such that $M(\omega)\subseteq M(\alpha'\alpha_1)\cup M(\mu\beta)$ and $(\mu\beta)^\omega=(\mu\beta)^{-1}$ and $\alpha_i^\omega=\alpha_i^{-1}$ for all $\alpha_i\in\{\alpha'\alpha_1\}$. As with Remark\ \ref{pan4-00004}, there exists a permutation $\gamma\in S_n$ such that $\alpha^\gamma=\alpha^{-1}$ and $\beta^\gamma=\beta^{-1}$. In particular, $\alpha_i^{\gamma}=\alpha_i^{-1}$ for each $\alpha_i\in\{\alpha'\}$.

Case 2: there exist two transitive cycle factors chains of $\beta$ on $\alpha$, and the second transitive cycle factors chain is $\alpha_{k+1}\rightarrow\alpha_{k+2}\rightarrow\cdot\cdot\cdot\rightarrow\alpha_l$. Regarding $[\alpha'\alpha_1,\mu\beta]$. Since $\mu\in C_{S_n}(\alpha_1\alpha_2\cdot\cdot\cdot\alpha_k)$ and $M(\mu)=M(\alpha_1\alpha_2\cdot\cdot\cdot\alpha_k)$, we see that $\alpha_{k+1}\rightarrow\alpha_{k+2}\rightarrow\cdot\cdot\cdot\rightarrow\alpha_l$ is the unique transitive cycle factors chain of $\mu\beta$ on $\alpha'\alpha_1$. According to Case 1 and Lemma\ \ref{pan2-8}, it follows that there exists a permutation $\omega'$ such that $M(\omega')\subseteq M(\alpha'\alpha_1)\cup M(\mu\beta)$, $(\mu\beta)^{\omega'}=(\mu\beta)^{-1}$ and $(\alpha'\alpha_1)^{\omega'}=(\alpha'\alpha_1)^{-1}$ and $\alpha_i^{\gamma}=\alpha_i^{-1}$ for each $\alpha_i\in(\{\alpha''\}\cup\{\alpha_1\})$ where $\{\alpha''\}=\{\alpha'\}\setminus\{\alpha_{k+1},\alpha_{k+2},...,\alpha_l\}$. Applying Remark\ \ref{pan4-00004}, there exists a permutation $\gamma\in S_n$ such that $\alpha^\gamma=\alpha^{-1}$ and $\beta^\gamma=\beta^{-1}$. In particular, $\alpha_i^{\gamma}=\alpha_i^{-1}$ for each $\alpha_i\in\{\alpha''\}$.

Case 3: there exist three transitive cycle factors chains of $\beta$ on $\alpha$, and the third transitive cycle factors chain is $\alpha_{l+1}\rightarrow\alpha_{l+2}\rightarrow\cdot\cdot\cdot\rightarrow\alpha_m$. Considering $[\alpha'\alpha_1,\mu\beta]$. In a similar manner, we know that there are two transitive cycle factors chains of $\mu\beta$ on $\alpha'\alpha_1$, those are $$\alpha_{k+1}\rightarrow\alpha_{k+2}\rightarrow\cdot\cdot\cdot\rightarrow\alpha_l~{\rm{and}}~\alpha_{l+1}\rightarrow\alpha_{l+2}\rightarrow\cdot\cdot\cdot\rightarrow\alpha_m.$$ According to Case 2 and Lemma\ \ref{pan2-8}, it follows that there exists a permutation $\omega''$ such that $M(\omega'')\subseteq M(\alpha'\alpha_1)\cup M(\mu\beta)$, $(\mu\beta)^{\omega''}=(\mu\beta)^{-1}$ and $\alpha'^{\omega''}=\alpha'^{-1}$ and $\alpha_1^{\omega''}=\alpha_1^{-1}$. Then by using Remark\ \ref{pan4-00004}, there exists a permutation $\gamma\in S_n$ such that $\alpha^\gamma=\alpha^{-1}$ and $\beta^\gamma=\beta^{-1}$. \qed

So far, we have confirmed Conjecture\ \ref{pan1-0} and thus we obtain the following theorem.

\begin{thm}\label{pan4-111}\normalfont
Let $\alpha$ and $\beta$ be two permutations in $S_n$. If the commutator $[\alpha,\beta]$ has at least $n-4$ fixed points, then there exists a permutation $\gamma\in S_n$ such that $\alpha^\gamma=\alpha^{-1}$ and $\beta^\gamma=\beta^{-1}$.
\end{thm}

Here, we give an example which shows that $k=4$ is the best bound in the sense of $|M([\alpha,\beta])|\leq k$ implies that there exists a permutation which simultaneously conjugates $\alpha$ and $\beta$ onto their respective inverses.
\begin{Exam}\label{pan8-18}\normalfont
Given $\alpha=(4,3,2,1)$ and $\beta=(3,2,1,5,4,6,7)$. Then we see $[\alpha,\beta]=(2,3,4,5,6)$. Noticing that there do not exist a permutation $\omega$ such that $\alpha^\omega=\alpha^{-1}$ and $\beta^\omega=\beta^{-1}$.
\end{Exam}

\section{Acknowledgement}

We are very grateful to the anonymous referees for their useful suggestions and comments for improvement. We also thank J. K$\rm{\ddot{o}}$nig and D. Neftin very much for their help.

\end{document}